%% file: paper_REGINA.tex
\documentclass[hidelinks,onefignum,onetabnum]{siamart251216}
\usepackage{multirow} 
\usepackage{comment}

\input{ex_shared}

\ifpdf
\hypersetup{
  pdftitle={AG-BFGS: Foundations for Accelerated Geometry Learning },
  pdfauthor={R. Sousa, D. Vargas, E. Wanner, and J. Knowles}
}
\fi

\begin{document}

\maketitle

\begin{abstract}
We introduce an accelerated geometry-learning framework for BFGS that incorporates momentum into the learning of Hessian inverse approximations. We prove that the resulting accelerated updates satisfy explicit stability conditions and preserve positive definiteness within admissible parameter regions. Numerical experiments demonstrate consistent improvements over the corresponding non-accelerated schemes while maintaining numerical stability. These results provide the first rigorous theoretical foundation for accelerated quasi-Newton learning and bridge the gap between empirical observations and mathematical guarantees for accelerated second-order methods.

\end{abstract}


\noindent

\input{introduction}

\input{accelerated}

\input{alternative_convergence}

\input{interpretation_param}
\input{experiments_new}
\input{conclusion}






\section*{Acknowledgments}

\bibliographystyle{siamplain}
\bibliography{references}
\end{document}

%% file: ex_shared.tex
\usepackage{lipsum}
\usepackage{amsfonts}
\usepackage{graphicx}
\usepackage{subcaption}
\usepackage{epstopdf}
\usepackage{algorithmic}
\usepackage{booktabs}

\ifpdf
  \DeclareGraphicsExtensions{.eps,.pdf,.png,.jpg}
\else
  \DeclareGraphicsExtensions{.eps}
\fi

\newsiamremark{remark}{Remark}
\newsiamremark{hypothesis}{Hypothesis}
\crefname{hypothesis}{Hypothesis}{Hypotheses}
\newsiamthm{claim}{Claim}
\newsiamremark{fact}{Fact}
\crefname{fact}{Fact}{Facts}

\headers{Accelerated Geometry BFGS}{R. Sousa, D. Vargas, E. Wanner, and J. Knowles}

\title{AG-BFGS: Foundations for Accelerated Geometry Learning\thanks{Submitted to arXiv 18 November, 2026
\funding{This work was partially funded by Fapemig, CAPES, and CNPq, Brazil.}}}

\author{Regina Carla Sousa\thanks{Dept of Mathematics, UFOP, Brazil 
  (\email{regina\_carla@ufop.edu.br} 
  ).}
\and Denis Vargas\thanks{PPGMMC, CEFET-MG, Brazil 
  (\email{denis.vargas@cefetmg.br}).}
\and Elizabeth Wanner \thanks{CSDT School, Aston University - PPGMMC, CEFET-MG  
  (\email{e.wanner1@aston.ac.uk}).}
  \and Joshua~Knowles \thanks{SLB Cambridge Research, High Cross, Madingley Road, Cambridgeshire, UK.  
  } }

\usepackage{amsopn}


%% file: introduction.tex
\section{Introduction} \label{sec:intro}

Quasi-Newton methods form a class of iterative algorithms for unconstrained smooth minimization that seek to combine the fast local convergence of Newton’s method with the lower computational cost of first-order approaches. Rather than explicitly computing the Hessian matrix (or its inverse, $H_k^{-1}$), these methods iteratively construct curvature approximations using gradient evaluations and differences between successive iterates. Among these methods, the BFGS algorithm  \cite{broyden1967,fletcher1970new,goldfarb1970family, shanno1970conditionin} has emerged as a particularly influential approach due to its robustness, strong practical performance, and superlinear local convergence behavior \cite{nocedal2006numerical}. Rather than computing exact Hessian matrices, BFGS recursively constructs curvature approximations using gradient information only, thereby avoiding the prohibitive computational cost associated with second-order derivatives in large-scale applications.

From a geometric perspective, BFGS may be interpreted as an iterative metric-learning process. At each iteration, the algorithm updates a symmetric positive definite approximation of the inverse Hessian matrix, which defines a local metric adapted to the geometry of the objective function. Consequently, the optimization process involves two coupled dynamical systems: the evolution of the iterates and the evolution of the underlying metric approximation. 
In high-dimensional settings, poorly conditioned landscapes, or multimodal functions, the performance of classical BFGS may degrade significantly, often leading to slow convergence or stagnation \cite{dai2002convergence,dai2013perfect,mascarenhas2004bfgs}. 
These limitations have motivated a broad range of extensions of quasi-Newton methods. Limited-memory formulations, such as L-BFGS, have been developed to make quasi-Newton techniques practical for large-scale optimization problems by storing only a small number of curvature pairs \cite{nocedal1980,nocedal1989,moritz2026,BYRD2019}. Richer curvature representations have also been proposed. In particular, Gao and Goldfarb \cite{Gao+2018} introduced Block BFGS, which incorporates multiple curvature directions simultaneously and achieves global and superlinear convergence under assumptions similar to those of classical BFGS. Beyond the Euclidean setting, Huang et al. \cite{Huang+18} developed a Riemannian BFGS method for nonconvex optimization on manifolds, demonstrating that quasi-Newton ideas can be successfully extended to optimization problems with non-Euclidean geometry. For large-scale machine learning applications, stochastic quasi-Newton methods have been proposed to cope with noisy gradient information and massive datasets \cite{gunter2007,byrd2016stochastic,wang2017stochastic,schraudolph2007stochastic,guo2023overview}; in particular, Wang et al. \cite{wang2017stochastic} established convergence guarantees for stochastic quasi-Newton schemes applied to nonconvex stochastic optimization.

Although these developments have substantially broadened the applicability of quasi-Newton methods, they primarily address improvements in curvature representation, extensions to non-Euclidean geometries, or robustness to stochastic information. The problem of accelerating the geometry-learning process itself, namely the evolution of the Hessian approximation, remains largely unexplored. While some accelerated quasi-Newton approaches have been proposed, a rigorous understanding of the resulting Hessian approximation dynamics, including stability, preservation of positive definiteness, and convergence guarantees, is still missing.

Acceleration techniques based on inertial dynamics and momentum mechanisms have become central in modern optimization. Beginning with Polyak’s heavy-ball method \cite{polyak1964some} and further developed through Nesterov’s accelerated gradient framework \cite{nesterov1983method}, momentum-based strategies have achieved significant improvements in convergence speed and transient optimization behavior. Their extension to quasi-Newton methods has attracted increasing attention \cite{su2016differential,gower2018accelerated,botev2017nesterov}. In particular, Gower et al. \cite{gower2018accelerated} introduced an accelerated BFGS framework based on inertial inverse-Hessian updates, providing an important step toward incorporating acceleration into quasi-Newton geometry learning.

Nevertheless, existing accelerated optimization methods predominantly introduce inertia into the evolution of the optimization variables, while the curvature approximation continues to evolve through standard recursive update mechanisms. Even in accelerated quasi-Newton schemes, Hessian approximations typically evolve through first-order update rules without explicit inertial memory in the metric dynamics.  The approach proposed in \cite{gower2018accelerated} provides empirical evidence that accelerated inverse-Hessian updates can improve performance over classical BFGS; however, a complete theoretical characterization of accelerated metric dynamics is still absent.

Motivated by these observations, we investigate the incorporation of inertia directly into the evolution of the quasi-Newton metric. From a dynamical perspective, the proposed approach introduces acceleration not only in the optimization trajectory but also in the evolution of the local geometry itself, leading naturally to the interpretation of the method as an inertial metric-learning process on the manifold of symmetric positive definite matrices. However, introducing inertia into quasi-Newton updates is mathematically nontrivial. The effectiveness of BFGS relies on preserving key structural properties, including symmetry, positive definiteness, and secant consistency. Naive inertial modifications may compromise these properties, potentially destroying descent directions and destabilizing the optimization process. Consequently, accelerated Hessian dynamics require careful analysis to ensure that the metric sequence remains well defined, stable, and numerically reliable. More broadly, despite the encouraging results reported in \cite{gower2018accelerated}, the theoretical understanding of deterministic accelerated quasi-Newton methods remains incomplete. In particular, convergence guarantees for deterministic accelerated BFGS schemes, the preservation of classical quasi-Newton properties such as positive definiteness and superlinear convergence, and the stability of inertial metric dynamics remain largely open research questions.

The main contributions of this work are summarized as follows. First, we propose an inertial BFGS update in which momentum is incorporated directly into the Hessian approximation dynamics, 
dubbed Accelerated Geometry BFGS (AG-BFGS). Second, we establish conditions under which the proposed update preserves symmetry and positive definiteness of the metric sequence. Third, we analyze the convergence properties of the resulting optimization method under standard smoothness assumptions and derive conditions ensuring stability of the inertial quasi-Newton dynamics. 
Finally, numerical experiments on 
several logistic regression problems demonstrate that the proposed accelerated metric dynamics can improve transient convergence behavior and enhance robustness in challenging optimization landscapes.

The remainder of this paper is organized as follows. Section \ref{sec:accelerated} presents the two main components of the proposed approach. Section \ref{sec:alter_convergence} analyzes the convergence properties of the proposed algorithm, establishing the conditions under which the method preserves the stability and convergence behavior of the classical BFGS method, as well as discussing parameter selection. Section \ref{sec:interpretation} analyze the acceleration parameters. Section \ref{sec:experiments} evaluates the performance of the proposed optimization algorithm on a set of logistic regression problems comparing with a classical BFGS and the approach proposed in \cite{gower2018accelerated}. Furthermore, the proposed accelerated approach is compared with its classical counterpart in a set of ill-conditioned and unimodal functions. Finally, Section \ref{sec:conclusion} concludes the paper and outlines the directions for future work.

%% file: accelerated.tex
\section{Accelerated Geometry BFGS}
\label{sec:accelerated}

The classical BFGS algorithm \cite{broyden1967,fletcher1970new,goldfarb1970family, shanno1970conditionin} iteratively builds an approximation $X_k \approx H_k^{-1}$ of the inverse Hessian matrix using only first-order information. The update is formulated to satisfy the secant condition, ensuring consistency between the displacement in the decision space and the observed variation in the gradient.
Consider a continuously differentiable objective function
$f:\mathbb{R}^n \rightarrow \mathbb{R}$,
whose gradient $\nabla f$ is $L$-Lipschitz continuous on a compact set containing all iterates $\{x_k\}$ generated by the algorithm.

Define the vectors
\begin{equation}
    s_k = x_{k+1} - x_k, 
\qquad
\zeta_k = \nabla f(x_{k+1}) - \nabla f(x_k),
\label{eq:s_zeta}
\end{equation}

\noindent in which $s_k$ denotes the step in the decision variable and $\zeta_k$ the corresponding change in the gradient. The secant condition imposes
\[
X_{k+1} \zeta_k = s_k,
\]
while symmetry, $X_{k+1} = X_{k+1}^\top$, is enforced to preserve the structure of the Hessian approximation.

The BFGS update can be characterized as the solution to the variational problem
\begin{equation}
X_{k+1} = \arg\min_X \|X - X_k\|_F^2 
\quad \text{subject to} \quad X \zeta_k = s_k, \; X = X^\top,
\label{eq:bfgs_opt}
\end{equation}
whose closed-form solution is given by
\begin{equation}
X_{k+1} 
= \frac{s_k s_k^\top}{s_k^\top \zeta_k}
+ \left(I - \frac{s_k \zeta_k^\top}{s_k^\top \zeta_k}\right)
X_k
\left(I - \frac{\zeta_k s_k^\top}{s_k^\top \zeta_k}\right).
\label{eq:bfgs_closed_form}
\end{equation}
The first term in Equation~\eqref{eq:bfgs_closed_form} enforces the secant condition, while the second provides the minimal Frobenius-norm correction, $\|A\|_F = \sqrt{\operatorname{Tr}(A^\top A)}$, preserving symmetry and, under the curvature condition $\zeta_k^\top s_k > 0$, the positive definiteness of $X_{k+1}$.

This formulation serves as the basis for the extensions developed in this work. In particular, the variational perspective in Equation~\eqref{eq:bfgs_opt} enables the incorporation of additional mechanisms directly into the update of $X_k$. In what follows, this structure is used to  introduce acceleration terms inspired by Nesterov’s method \cite{nesterov1983method}
 while preserving the favorable convergence properties of the quasi-Newton framework.

Acceleration methods based on inertial dynamics, proposed by \cite{nesterov1983method}, augment the standard BFGS quasi‑Newton framework with momentum-inspired acceleration to improve convergence speed and reduce stagnation, particularly in ill‑conditioned or weakly convex problems. In their classical form, 
\begin{equation}
\begin{array}{l}
y_k = x_k + \rho (x_k - x_{k-1}) \\
x_{k+1} = y_k - \alpha \nabla f(y_k),
\end{array}
\label{eq:clasnesterov}
\end{equation}
these accelerated methods introduce a momentum step ($\rho$) that looks ahead before evaluating gradients or updating the inverse Hessian. For analytical convenience and improved numerical stability we introduce a third sequence ($v_k$) that captures a smoothed inertial state by aggregating information from previous iterates. Defining $v_0 = x_0$ and for $k \geq 0$ set
$v_{k+1} = (1 - \beta) v_k + \beta x_{k+1}$,
with weights $\beta \in (0,1]$, Equations (\ref{eq:clasnesterov}) can be written as
\begin{equation}
\begin{array}{l}
y_k = (1 - \rho) x_k + \rho (v_k) \\
x_{k+1} = y_k - \alpha \nabla f(y_k) \\
v_{k+1} = (1 - \beta)v_k + \beta x_{k+1}.
\end{array}
\label{eq:nesterov}
\end{equation}
The parameter $\rho$ controls the extrapolation between the current state $x_k$ and the inertial state $v_k$, while $\alpha$ denotes the step size, and $\beta$ regulates the update of the dynamic memory $v_k$. 
Note that the gradient correction term is weighted by the product $\beta\alpha$, implying that the influence of the inertial memory and the magnitude of the descent correction are intrinsically coupled through the parameter $\beta$.

In this work, rather than applying acceleration directly to the decision variables, we adopt a different perspective by transferring the principles of acceleration to the space of matrices that approximate the inverse Hessian. Specifically, inertial dynamics is incorporated into the sequence ${X_k}$, where $X_k \approx H_k^{-1}$, following the approach of \cite{gower2018accelerated}. In this setting, the inertial component is no longer a vector but an auxiliary matrix that accumulates historical information from past approximations of $H^{-1}$. Consequently, acceleration operates directly on the sequence ${X_k}$, enhancing the evolution of inverse Hessian approximations rather than acting on the decision vector $x_k$.

Adapting Nesterov’s acceleration to the matrix setting requires conceptual and algebraic adjustments. In contrast to vector-valued variables, inverse Hessian approximations lie in the cone of symmetric positive definite matrices, denoted by $\mathbb{S}{++}^n$. This set is not a vector space since arbitrary linear combinations of SPD matrices may produce indefinite matrices. As a result, the standard linear extrapolation used in the classical method 
\begin{equation*}X_k + \rho (X_k - X_{k-1}),
\end{equation*} 
in general does not preserve positive definiteness and cannot be applied directly to update the sequence ${X_k}$.

The aim remains to introduce an inertial mechanism that anticipates the evolution of the quasi-Newton iterates. However, the nonlinear geometry of $\mathbb{S}_{++}^n$ prevents the use of a single momentum term as in the classical Nesterov scheme. Instead, the matrix formulation employs two complementary parameters: $\rho \in (0,1)$, which governs a stable interpolation between successive inverse Hessian approximations, and $\beta \in [0,1)$, which controls inertial smoothing in the auxiliary sequence. This separation ensures that positive definiteness is maintained while mitigating the risk of numerical instability due to overly aggressive extrapolation.

Unlike Nesterov’s accelerated method in the iterated dynamics, the proposed accelerated quasi-Newton formulation does not explicitly employ a step-size parameter associated with the gradient, such as the coefficient $\alpha$ in  $x_{k+1} = y_k - \alpha \nabla f(y_k)$.
This is because, in the quasi-Newton setting, the matrix $X_k$
itself acts as an adaptive scaling and preconditioning mechanism for the descent direction, implicitly incorporating local curvature information of the objective function. As a result, the effective magnitude of the update is no longer governed by an explicit gradient step-size parameter, but instead emerges from the approximation of the inverse Hessian together with the line-search procedure used to update $x_k$.

The auxiliary sequence ${V_k}$ accumulates, in a smoothed fashion, the historical information from the update process. At each iteration, $V_k$ and $X_k$ are combined via a convex combination, ensuring that the intermediate point remains within the cone $\mathbb{S}{++}^n$:
\begin{equation}
Y_k = (1-\rho)X_k + \rho V_k,
\label{eq:yk_iter}
\end{equation}
in which $\rho$ controls the inertial contribution and plays a role analogous to the momentum parameter in the vector setting, while being restricted to convex interpolation to guaranty that $Y_k \in \mathbb{S}{++}^n$, i.e. $Y_k \succ 0$. In this way, $Y_k$ serves as a matrix analog of Nesterov’s extrapolated point $y_k$, anticipating the direction of motion in the metric space without compromising the geometric consistency of the method.

The quasi-Newton update is then performed at $Y_k$, producing a new approximation $X_{k+1}$:
\begin{equation}
X_{k+1} = \Phi_{\text{BFGS}}(Y_k;, s_k, \zeta_k),
\label{eq:xk_iter}
\end{equation}
in which $\Phi_{\text{BFGS}}$ denotes the BFGS transformation 
\begin{equation}
\Phi_{\text{BFGS}}(Y_k;, s_k, \zeta_k)
= \frac{s_k s_k^\top}{s_k^\top \zeta_k}
+ \left(I - \frac{s_k \zeta_k^\top}{s_k^\top \zeta_k}\right)
Y_k
\left(I - \frac{\zeta_k s_k^\top}{s_k^\top \zeta_k}\right).
\label{eq:bfgs_closed_form1}
\end{equation}
This step preserves positive definiteness whenever $s_k^\top \zeta_k > 0$, a standard assumption in quasi-Newton methods. Consequently, the update ensures that $X_{k+1} \in \mathbb{S}_{++}^n$, even if $Y_k$ is close to the boundary of the cone.

Adapting the accelerated structure to the quasi-Newton setting motivates the introduction of an additional parameter $\gamma$, responsible for decoupling the inertial memory from the correction induced by the BFGS update. The acceleration BFGS is given by
\begin{equation}
\begin{array}{l}
x_{k+1} = x_k - \eta_k X_k \nabla{f} (x_k), \\
Y_k = (1-\rho)X_k + \rho V_k, \\
X_{k+1} = \Phi_{\text{BFGS}}(Y_k;, s_k, \zeta_k), \\
V_{k+1} = (1-\beta)V_k + (\beta-\gamma)Y_k + \gamma X_{k+1}.
\end{array}
\label{eq:acc_bfgs}
\end{equation}
The parameter $\eta_k$ is the step-size, while the parameter $\rho \in [0,1)$ controls the matrix extrapolation used in the construction of $Y_k$. The parameter $\beta$ controls the degree of inertial damping, while the parameter $\gamma$ governs the extent to which the correction induced by the BFGS transformation is incorporated into the accelerated dynamics. The appropriate ranges for $\beta$ and $\gamma$ will be derived in Section \ref{sec:alter_convergence}.

Equations (\ref{eq:acc_bfgs}) define a discrete and numerically stable acceleration scheme in which:
\begin{enumerate}
\item \label{it:item1} \textbf{Convex combination} of $X_k$ and $V_k$ preserves positive definiteness, ensuring $Y_k \in \mathbb{S}{++}^n$;
\item \label{it:item2} \textbf{The quasi-Newton update} $\Phi{\text{BFGS}}$ maintains symmetry and positive definiteness of $X_{k+1}$ under the condition $s_k^\top \zeta_k > 0$;
\item \label{it:item3} \textbf{The inertial update} of $V_k$ smooths the dynamics, controlling the propagation of inertia and preventing oscillations in ill-conditioned regions.
\end{enumerate}

From a conceptual standpoint, the accelerated quasi-Newton iteration can be viewed as a two-layer mechanism: the inner layer, defined by $\Phi_{\text{BFGS}}$, enforces geometric consistency within the SPD cone, while the outer layer, governed by the updates of $Y_k$ and $V_k$, introduces a Nesterov-type inertial effect in a fully discrete and matrix-compatible form. This yields an iterative scheme that combines the adaptability of quasi-Newton methods with controlled acceleration, without relying on continuous-time analysis or parameters derived from differential equations. 

Algorithm~\ref{alg:bfgs_ac_stoc} summarizes the proposed formulation dubbed Accelerated Geometry BFGS (AG-BFGS).  Parameters $\rho$ and $\beta$ are chosen within the ranges 
established by the stability analysis, ensuring preservation of structural properties in the cone $\mathbb{S}_{++}^n$.

\begin{algorithm}[H]
\caption{Accelerated Geometry BFGS  (AG-BFGS)}
\label{alg:bfgs_ac_stoc}
\begin{algorithmic}[1]
\REQUIRE $x_0$, $X_0 = I$, $U_0 = I$, 
$\eta_k$, $\rho$, $\beta$, $\gamma$
\STATE \textbf{Evaluate} $f(x_0)$ and $g_0 = \nabla f(x_0)$
\FOR{$k = 0,1,2,\dots$ \textbf{until convergence}}
    \STATE $d_k \gets -X_k g_k$ 
    \STATE \textbf{Update iterate}: $x_{k+1} = x_k + \eta_k d_k$
    \STATE \textbf{Compute} $g_{k+1} = \nabla f(x_{k+1})$
    \STATE \textbf{Set} $s_k = x_{k+1} - x_k$, \quad $y_k = g_{k+1} - g_k$
    \IF{$y_k^\top s_k > 0$}
        \STATE $Y_k \gets (1-\rho)X_k + \rho V_k$
        \STATE $\phi_k \gets 1/(y_k^\top s_k)$
        \STATE $W_k \gets I - \phi_k s_k y_k^\top$
        \STATE $X_{k+1} \gets W_k Y_k W_k^\top + \phi_k\, s_k s_k^\top$
        \STATE $V_{k+1} \gets (1-\beta)V_k + (\beta -\gamma)Y_k + \gamma X_{k+1}$
    \ELSE
        \STATE $X_{k+1} \gets X_k$, \quad $V_{k+1} \gets V_k$
    \ENDIF
\ENDFOR
\end{algorithmic}
\end{algorithm}

%% file: alternative_convergence.tex
\section{Convergence Analysis}
\label{sec:alter_convergence}

In this section, we analyze the convergence properties of the proposed Accelerated Geometry BFGS (AG-BFGS) algorithm. The analysis establishes conditions under which the method retains the stability and convergence guarantees of classical BFGS, while incorporating inertial momentum-like effects in the Hessian updates. In particular, we aim to characterize admissible regions for $(\rho,\beta,\gamma)$ that ensure both the dynamical stability of the accelerated system and the preservation of the structural properties of the quasi-Newton approximations.

\subsection{Iterative Analysis and Preservation of Structural Properties}
\label{subsec:analise_iterativa}

To establish the structural and stability properties of the proposed Accelerated Geometry BFGS (AG-BFGS) scheme described in Algorithm \ref{alg:bfgs_ac_stoc} (Section \ref{sec:accelerated}), we present a set of auxiliary results. These results provide sufficient conditions for the preservation of positive definiteness and for the monotonic decrease of a Lyapunov functional, $\mathcal{E}_k$, associated with the sequences $\{X_k, V_k, Y_k, x_k\}$.

A central requirement is the boundedness of the eigenvalues of the metric sequence $X_k,~ Y_k,~V_k$. In particular, it is necessary to establish the existence of constants $0< m\leq M$ such that  $mI \preceq X_k, Y_k, V_k \preceq MI$ for some $0 < m \le M$. These bounds ensure that the metric approximation remains uniformly positive definite and well conditioned, preventing instability of the search directions and potential divergence caused by the inertial Hessian dynamics.

\begin{lemma}[Boundedness of Eingenvalues]
\label{lem:convex_spd_bounds}
Let $X_k, V_k \in \mathbb{S}_{++}^n$ satisfy $mI \preceq X_k, V_k \preceq MI$ for some $0 < m \le M$. Then, for any $\rho \in (0,1)$,
\begin{equation*}
Y := (1-\rho)X_k + \rho V_k \in \mathbb{S}_{++}^n
\quad\text{and}\quad
mI \preceq Y \preceq MI.
\end{equation*}

\begin{proof}
For any $z \neq 0$, we have that
\begin{equation*}
z^\top Y z = (1-\rho) z^\top X_k z + \rho z^\top V_k z > 0.
\end{equation*}
Moreover,
\begin{equation*}
\lambda_{\min}(Y) \ge (1-\rho)\lambda_{\min}(X_k) + \rho\lambda_{\min}(V_k) \ge m,
\end{equation*}
and
\begin{equation*}
\lambda_{\max}(Y) \le (1-\rho)\lambda_{\max}(X_k) + \rho\lambda_{\max}(V_k) \le M,
\end{equation*}
which concludes the proof.
\end{proof}
\end{lemma}

An entirely analogous argument applies to the auxiliary variable $V_k$, whose update is also given by a convex combination of matrices in $\mathbb{S}_{++}^n$. Consequently, both $Y_k$ and $V_{k+1}$ preserve positive definiteness by construction, ensuring that all matrix variables remain in the cone $\mathbb{S}_{++}^n$ throughout the iterations.

Under the curvature condition $s_k^\top \zeta_k > 0$, the BFGS transformation preserves positive definiteness, thereby completing the structural stability cycle between extrapolation and quasi-Newton update.

\begin{lemma}[Positive Definiteness]
\label{lem:bfgs_spd}
Let $Y_k \in \mathbb{S}_{++}^n$ and $s, \zeta \in \mathbb{R}^n$ be such that $\zeta^\top s > 0$. The classical BFGS update
\begin{equation*}
X_{k+1}
= \Big(I-\tfrac{s \zeta^\top}{\zeta^\top s}\Big)\,Y_k\,\Big(I-\tfrac{\zeta s^\top}{\zeta^\top s}\Big)
  + \tfrac{s s^\top}{\zeta^\top s}
\end{equation*}
satisfies $X_{k+1} \in \mathbb{S}_{++}^n$.

\begin{proof}
Let $\phi := \frac{1}{\zeta^\top s} > 0$ and define $W := I - \phi\,\zeta s^\top$. 

Then $W^\top = I - \phi\, s\zeta^\top$, and the update can be written as
\begin{equation} \begin{array}{lcl}
X_{k+1} &
=& \big(I-\phi\, s\zeta^\top\big)\,Y_k\,\big(I-\phi\,\zeta s^\top\big) + \phi\,ss^\top  \\
&=& W^\top Y_k W + \phi\,ss^\top .
\end{array}
\nonumber
\end{equation}
For any $z \neq 0$, we have that
\begin{equation*}
z^\top X_{k+1} z = (Wz)^\top Y_k (Wz) + \phi (s^\top z)^2.
\end{equation*}
If $Wz \neq 0$, the first term is strictly positive since $Y_k \succ 0$. If $Wz = 0$, then $z = \phi \zeta (s^\top z)$, which implies $s^\top z \neq 0$, and hence the second term is strictly positive. Therefore, $z^\top X_{k+1} z > 0$ for all $z \neq 0$, proving that $X_{k+1} \in \mathbb{S}_{++}^n$. Symmetry follows from the symmetry of $Y_k$ and the update structure.
\end{proof}
\end{lemma}

At this point, the accelerated scheme \eqref{eq:acc_bfgs}, together with the BFGS transformation \eqref{eq:bfgs_closed_form1}, preserves the structural properties listed in Section~\ref{sec:accelerated}: the convex combination of $X_k$ and $V_k$ remains positive definite so $Y_k \in \mathbb{S}_{++}^n$ , and the quasi‑Newton update $\Phi{\text{BFGS}}$ yields a symmetric, positive and positive-definite $X_{k+1}$ provided $s_k^\top \zeta_k > 0$. Therefore, throughout the remainder of the analysis, we assume that the curvature condition $s_k^\top \zeta_k >0$,
holds for all iterations and that there exist constants $0<m\le M<\infty$ such that
\begin{equation}
mI \preceq X_k,V_k,Y_k \preceq MI,
\label{eq:auto_val}
\end{equation}
for every $k\ge0$.

Furthermore, since \(f:\mathbb{R}^n\rightarrow\mathbb{R}\) is continuously differentiable with \(L\)-Lipschitz continuous gradient on a compact set containing the iterates \(\{x_k\}\), we assume that the step sizes \(\eta_k\) are chosen so that the sufficient decrease condition
\begin{equation}
f(x_{k+1})
\le
f(x_k)-c_1\|\nabla f(x_k)\|^2,
\label{eq:suff_decrease}
\end{equation}
holds for some constant \(c_1>0\), as ensured, for example, by an Armijo line search or by sufficiently small step sizes.

Under these assumptions, the next lemma quantifies the discrepancy introduced by replacing the extrapolated matrix $Y_k$ with the corrected BFGS approximation $X_{k+1}$. This estimate plays a key role in the Lyapunov analysis, allowing the inertial component of the accelerated dynamics to be controlled and ultimately leading to the stability result established in Theorem~\ref{thm:Lyapunov_revised}.

\begin{lemma}[BFGS perturbation estimate]\label{lem:BFGS_perturb}
Under the assumptions \eqref{eq:auto_val} and \eqref{eq:suff_decrease}, there exists $A_1>0$ 
(depending only on $L,m,M$) 
such that for all $k$ 
\[\|X_{k+1}-Y_k\|_F^2 \le A_1 \eta_k^4 \|\nabla f(x_k)\|^4.\]
\end{lemma}

\begin{proof}[Proof Lemma]
By the quasi-Newton lemma (based on the BFGS update applied to $Y_k$), there exists a constant $A>0$, depending only on $m$, $M$, and $L$, such that
\begin{equation}\label{eq:bound-Xdiff-1}
\|X_{k+1}-Y_k\|_F
\;\le\;
A\,\|s_k\|\,\|\zeta_k\|.
\end{equation}

\noindent Using \eqref{eq:s_zeta} together with the Lipschitz continuity of the gradient, we obtain
\[
\|X_{k+1}-Y_k\|_F
\le
A\,\|s_k\|\,\|\zeta_k\|
\le
ALM^2\,\eta_k^2\,\|\nabla f(x_k)\|^2.
\]

\noindent Consequently, there exists a constant $A_1>0$ such that
\begin{equation}\label{eq:bound-Xdiff-2}
\|X_{k+1}-Y_k\|_F^2
\;\le\;
A_1\,\eta_k^4\,\|\nabla f(x_k)\|^4.
\end{equation}

\end{proof}


The following theorem provides the main stability result for the proposed AG-BFGS method. It shows that, under the smoothness assumptions stated above and for sufficiently small step sizes, there exists a Lyapunov function
$\mathcal{E}_k$
that decreases monotonically along the iterations. This function combines the objective value with a consistency measure between the matrix sequences $X_k$ and $V_k$, thereby capturing both optimization progress and stability of the accelerated dynamics.

The result further establishes that the monotone decrease of $\mathcal{E}_k$ is guaranteed whenever the acceleration parameters satisfy suitable admissibility conditions. Consequently, the theorem provides explicit stability regions for the parameters $\rho$, $\beta$, and $\gamma$, ensuring that the inertial mechanism remains compatible with the quasi-Newton corrections and does not destabilize the matrix dynamics.

\begin{theorem}[Numerical stability]
\label{thm:Lyapunov_revised}
Consider the accelerated scheme \eqref{eq:acc_bfgs} with the BFGS update \eqref{eq:bfgs_closed_form1}. Assume that the hypotheses stated above hold, in particular \eqref{eq:auto_val} and \eqref{eq:suff_decrease}. Let the parameters satisfy $0<\rho<1$, $0<\gamma<1$, and
\begin{equation}
\label{eq:beta_rho_cond}
2(1-\beta)^2(1-\rho)^2<1.
\end{equation}
For $\kappa>0$, define the Lyapunov function
\[
\mathcal{E}_k
:=
f(x_k)-f^\star
+
\kappa\|X_k-V_k\|_F^2.
\]

Then there exist constants $\kappa^\star>0$ and $c_2>0$ such that, for every $\kappa\in(0,\kappa^\star)$ and all $k$,
\begin{equation}
\label{eq:lyap_decrease}
\mathcal{E}_{k+1}
\le
\mathcal{E}_k
-
\frac{c_1}{2}\|\nabla f(x_k)\|^2
-
c_2\|X_k-V_k\|_F^2.
\end{equation}
Consequently, the sequence $\{\mathcal{E}_k\}$ is nonincreasing and $\|X_k-V_k\|_F \to 0$ as $ k\to\infty$.
\end{theorem}

\begin{proof}
Write $D_k:=X_k-V_k$ and $R_k:=X_{k+1}-Y_k$.

\noindent Step A - decomposition. 

\noindent From the definition,
\[
\mathcal{E}_{k+1}-\mathcal{E}_k = \underbrace{f(x_{k+1})-f(x_k)}_{\text{(I)}} \;+\; \kappa\big(\|X_{k+1}-V_{k+1}\|_F^2-\|D_k\|_F^2\big).
\]
By \eqref{eq:suff_decrease}, (I) $\le -c_1\|\nabla f(x_k)\|^2$. It remains to bound
\[
\Delta_k^{\mathrm{mat}}:=\|X_{k+1}-V_{k+1}\|_F^2-\|D_k\|_F^2.
\]

\noindent Step B - algebraic identity for $X_{k+1}-V_{k+1}$. 

\noindent Using the $V$-update and $Y_k=(1-\rho)X_k+\rho V_k$,
\begin{align}
X_{k+1}-V_{k+1}
&= (1-\beta)(X_{k+1}-V_k) + (\gamma-\beta)(Y_k-X_{k+1})\nonumber\\
&= (1-\beta)\big(R_k+(1-\rho)D_k\big) + (1-\gamma)R_k\nonumber\\
&= (1-\beta)(1-\rho)D_k + (1-\gamma)R_k. \label{eq:decomp_XminusV}
\end{align}

\noindent Step C - quadratic bound (Young inequality). 

\noindent For any matrices $A,B$ and any $\alpha>0$,
\[
\|A+B\|_F^2 \le (1+\alpha)\|A\|_F^2 + \Big(1+\frac{1}{\alpha}\Big)\|B\|_F^2.
\]
Take $A=(1-\beta)(1-\rho)D_k$, $B=(1-\gamma)R_k$, and $\alpha=1$ to get
\begin{equation}\label{eq:young_bound}
\|X_{k+1}-V_{k+1}\|_F^2 \le 2(1-\beta)^2(1-\rho)^2\|D_k\|_F^2 + 2(1-\gamma)^2\|R_k\|_F^2.
\end{equation}
Hence
\[
\Delta_k^{\mathrm{mat}} \le \big[2(1-\beta)^2(1-\rho)^2-1\big]\|D_k\|_F^2 + 2(1-\gamma)^2\|R_k\|_F^2.
\]
By condition \eqref{eq:beta_rho_cond} the coefficient $\omega:=2(1-\beta)^2(1-\rho)^2-1$ is negative; 
so that
\begin{equation}\label{eq:delta_preR}
\Delta_k^{\mathrm{mat}} \le -\omega\|D_k\|_F^2 + 2(1-\gamma)^2\|R_k\|_F^2.
\end{equation}


\noindent It follows from Lemma~\ref{lem:BFGS_perturb} that
\begin{equation}\label{eq:Rk_bound}
\|R_k\|_F^2 = \|X_{k+1}-Y_k\|_F^2 \le A_1\,\eta_k^4\,\|\nabla f(x_k)\|^4.
\end{equation}
Substituting \eqref{eq:Rk_bound} into \eqref{eq:delta_preR} gives
\[
\Delta_k^{\mathrm{mat}} \le -\omega\|D_k\|_F^2 + 2A_1(1-\gamma)^2\eta_k^4\|\nabla f(x_k)\|^4.
\]

\noindent Step D  - convert quartic to quadratic and choose $\kappa$. 

\noindent The compactness of the iterate sequence and the continuity of $\nabla f$ imply the existence of a uniform bound $G>0$ such that
\(
\|\nabla f(x_k)\|\le G, \forall k,
\)
as well as a constant $\eta_{\max}>0$ for which the following estimates hold whenever $\eta_k\le\eta_{\max}$. Hence
\[
\eta_k^4\|\nabla f(x_k)\|^4 \le \eta_{\max}^4 G^2 \|\nabla f(x_k)\|^2.
\]
Define $C := 2A_1(1-\gamma)^2\eta_{\max}^4 G^2$. Then
\[
\Delta_k^{\mathrm{mat}} \le -\omega\|D_k\|_F^2 + C \|\nabla f(x_k)\|^2.
\]
Returning to the decomposition (Step A) and using the sufficient decrease,
\begin{align*}
\mathcal{E}_{k+1}-\mathcal{E}_k
&\le -c_1\|\nabla f(x_k)\|^2 + \kappa\big(-\omega\|D_k\|_F^2 + C\|\nabla f(x_k)\|^2\big)\\
&= -\Big(c_1-\kappa C\Big)\|\nabla f(x_k)\|^2 - \kappa\omega\|D_k\|_F^2.
\end{align*}
Choose $0<\kappa\le \kappa^\star:=c_1/(2C_*)$ to ensure $c_1-\kappa C\ge c_1/2$. Then \eqref{eq:lyap_decrease} holds with $c_2:=\kappa\omega>0$.

\noindent Step E - consequences. 

\noindent The one-step decay \eqref{eq:lyap_decrease} implies $\{\mathcal{E}_k\}$ is nonincreasing and (since the second term is summable) $\|X_k-V_k\|_F\to0$. Summability of $\|\nabla f(x_k)\|^2$ also follows by summing \eqref{eq:lyap_decrease}.
\end{proof}

%% file: interpretation_param.tex
\section{Interpretation of the Acceleration Parameters}
\label{sec:interpretation}

An interesting feature of the proposed formulation is that the stability condition obtained from the Lyapunov analysis does not directly constrain the value of the parameter $\beta$ itself, but rather the effective magnitude of the inertia introduced into the system. Indeed, the dissipativity condition derived above (Equation (\ref{eq:beta_rho_cond})) can be written as



\begin{equation}
|1-\beta|
<
\frac{1}{\sqrt{2}(1-\rho)}.
\label{eq:cond_beta}
\end{equation}



From this perspective, the parameter $\rho$ acts as a moderator of the influence of $\beta$. When $\rho$ is close to zero, the factor $(1-\rho)$ remains close to one, and stability imposes a relatively narrow range for $\beta$. In this case,

\[
\rho\rightarrow0
\quad\Longrightarrow\quad
0.2929<\beta<1.7071.
\]

On the other hand, as $\rho$ approaches $1$, the factor $(1-\rho)$ tends to zero, significantly reducing the dynamical impact of $\beta$. As a consequence, the admissible region for $\beta$ becomes progressively larger,

\[
\rho\rightarrow1
\quad\Longrightarrow\quad
\beta\in\mathbb{R}.
\]

This result indicates that large values of $\beta$, including values greater than $1$, may be employed without compromising the stability of the accelerated system, provided that they are accompanied by sufficiently large values of $\rho$. In other words, the parameter $\rho$ acts as a damping mechanism capable of compensating for the increased inertia introduced by $\beta$.

Beyond dynamical stability, it is also necessary to preserve the geometric properties of the auxiliary sequence $V_k$. Expanding the update

\begin{equation}
V_{k+1}=(1-\beta)V_k+\beta Y_k-\gamma(Y_k-X_{k+1}),
\label{eq:Vk_expand}
\end{equation}
and substituting

\[
Y_k=(1-\rho)X_k+\rho V_k,
\]
yields

\begin{equation}
V_{k+1}
=
\bigl[1-\beta(1-\rho)-\gamma\rho\bigr]V_k
+
(\beta-\gamma)(1-\rho)X_k
+
\gamma X_{k+1}.
\label{eq:Vk_convexa}
\end{equation}

The coefficients in \eqref{eq:Vk_convexa} sum to one, so that $V_{k+1}$ can be interpreted as an affine combination of $V_k$, $X_k$, and $X_{k+1}$. To ensure that this combination is convex, and therefore automatically preserves positive definiteness whenever $V_k$, $X_k$, and $X_{k+1}$ belong to the cone $\mathbb{S}_{++}^n$, it suffices to impose

\begin{equation}
0\leq\gamma\leq\beta,
\qquad
\beta(1-\rho)+\gamma\rho\leq1.
\label{eq:cond_convexa}
\end{equation}

\noindent The second inequality implies

\begin{equation}
\beta
\leq
\frac{1-\gamma\rho}{1-\rho}.
\label{eq:beta_spd}
\end{equation}

Once again, the fundamental role of $\rho$ becomes apparent. When $\rho=0$, the classical restriction $\beta\leq1$ is recovered. However, for $\rho>0$, the upper bound increases progressively. For example, when $\rho=0.5$, one obtains approximately $\beta\leq2$, whereas for $\rho=0.9$, values as large as $\beta\leq10$ become admissible without violating the convexity of the update.

Therefore, the parameter $\rho$ plays a dual role in the proposed accelerated dynamics. From the stability viewpoint, it reduces the effective inertia associated with $\beta$. From the geometric viewpoint, it enlarges the admissible range of $\beta$ while preserving the membership of the sequence $V_k$ in the cone of symmetric positive-definite matrices. This interplay makes it possible to explore significantly more aggressive acceleration regimes than those typically considered in conventional momentum-based schemes, while simultaneously maintaining stability and geometric structure.

Theorem~\ref{thm:Lyapunov_revised} establishes Lyapunov stability for the proposed  accelerated geometry scheme. The decay condition, given by Equation \eqref{eq:lyap_decrease}, ensures that the energy sequence $\{\mathcal{E}_k\}$ decreases monotonically, jointly controlling the evolution of the primal variable $x_k$ and the matrix sequences $X_k$ and $V_k$. The dissipative term $\|X_{k+1}-V_k\|_F^2$ enforces alignment between the quasi-Newton update and the inertial state, implying that $X_k$ and $V_k$ asymptotically coincide. This mechanism acts as a damping effect, preventing oscillations induced by the inertial term. As a consequence, the method preserves positive definiteness, maintains uniform spectral bounds, and converges towards a stationary regime characterized by $\nabla f(x_k)\to 0$ and $X_k \approx V_k$, without compromising the structural properties of the BFGS update.

%% file: experiments_new.tex
\section{Numerical Experiments}
\label{sec:experiments}

In this section, we evaluate the performance of the proposed optimization algorithm, which extends the classical BFGS method by incorporating an accelerated geometry component. To this end, we considere three algorithmic variants in the computational study. The first corresponds to the standard BFGS algorithm, which serves as a baseline reference. The second is the accelerated BFGS method introduced in \cite{gower2018accelerated} (aBFGS). The third is the proposed AG-BFGS variant, which integrates the same accelerated geometry framework but is complemented by a convergence analysis that establishes admissible parameter regimes.

Although aBFGS and AG-BFGS are algorithmically equivalent, they differ in how the parameters are specified and interpreted. In particular, aBFGS is formulated in terms of the parameters $\mu$ and $\nu$, which in the original presentation are chosen ad-hoc without a formal characterization of their admissible range. In contrast, AG-BFGS is expressed in terms of parameters $\rho$, $\beta$, and $\gamma$, for which convergence conditions are explicitly derived, providing a theoretically justified parameter space.

The relationship between both parameterizations is given in \cite{gower2018accelerated} and using a straightforward transformation it follows that 
\begin{equation}
\mu = \frac{\beta}{\gamma},
\qquad
\nu = \frac{1}{\beta\gamma},
\qquad
\rho = \frac{\beta}{1+\beta}.
\label{eq:inverse_parameters}
\end{equation}
Moreover, the convergence analysis in Section~\ref{sec:alter_convergence} imposes the conditions
$0 < \gamma \leq \beta$
and
$\beta \leq \dfrac{1-\gamma\rho}{1-\rho}$,
which allow us to derive explicit bounds for the equivalent parameters. In particular, $\mu$ satisfies
\begin{equation}
1 \leq \mu \leq
\frac{1-\gamma\rho}{\gamma(1-\rho)},
\label{eq:mu_bound}
\end{equation}
while $\nu$ is bounded by
\begin{equation}
\frac{1-\rho}
{\gamma(1-\gamma\rho)}
\leq
\nu
\leq
\frac{1}{\gamma^2}.
\label{eq:nu_bound}
\end{equation}

Although relations given by Equation \eqref{eq:inverse_parameters} allow us to express the
aBFGS parameters $(\mu,\nu)$ in terms of $(\rho,\beta,\gamma)$, the bounds given 
in Equation \eqref{eq:mu_bound} and Equation \eqref{eq:nu_bound} should not be interpreted as
necessary conditions for the original aBFGS parameters. These bounds are
derived from the sufficient stability condition adopted in the proposed
AG-BFGS formulation. Therefore, parameter values previously suggested for
aBFGS may lie outside these intervals without implying inconsistency with
the original method. In the computational experiments, this relationship is used to ensure a consistent comparison between parameter settings while guaranteeing that all tested settings remain within theoretically admissible regimes. 

\subsection{Logistic Regression Problems}
\label{subsec:logistic}

We tested the three algorithms on four logistic regression problems using data from LIBSVM \cite{LIBSVM} (phishing, mushrooms, autralian, and splice datasets). To ensure comparability, we followed the experimental protocol provided in \cite{gower2018accelerated}. In all tests, the data was centered and normalized, including a bias term. The regularization parameter was set to $\lambda = \frac{1}{m}$, in which $m$ is the number of data points. The step-size, $\eta$, was chosen and kept fixed throughout the search. The stopping criteria was set to a residual of $1.00E-16$ or reaching the dataset-specific maximum number of iterations.

The aBFGS parameters, $\mu$ and $\nu$, used for each dataset were the values recommended in \cite{gower2018accelerated}, since the authors report that these parameters were manually selected to produce the best results (see Table \ref{tab:parameters_values} ). On the other hand, for AG-BFGS, $\rho$, $\beta$, and $\gamma$, as defined in Section~\ref{sec:alter_convergence}, were chosen using a generic, problem-independent strategy. Specifically, $\rho$ was set to the midpoint of its admissible interval. Given this value of $\rho$, $\beta$ was then set to the midpoint of its corresponding admissible interval, and, given the resulting values of $\rho$ and $\beta$, $\gamma$ was set to the midpoint of its corresponding admissible interval. This procedure resulted in the parameter values $\rho = [\text{0.35}]$, $\beta = [\text{0.75}]$, and $\gamma = [\text{0.5}]$, which were used for all datasets. This strategy requires no prior information about the problems and avoids any computational effort associated with parameter tuning or exhaustive search over the admissible parameter ranges.


For both methods, we employed the fixed step-size strategy as proposed in \cite{gower2018accelerated}, 
$\eta \in (0.50,~0.80,~1.00,~0.55)$.
The maximum number of iterations ($MaxIter$) was set to  $(120,~160,~250,~80)$ for the Australian, Phishing, Mushroom, and Splice datasets, respectively. 

\begin{table}[htbp]
\centering
\caption{Parameter values for aBFGS as proposed in \cite{gower2018accelerated}.}
\label{tab:parameters_values}
\begin{tabular}{lcc}
\toprule
\multirow{2}{*}{Dataset}  
& \multicolumn{2}{c}{aBFGS}  \\
\cmidrule(lr){2-3}
& $\mu$ & $\nu$ \\
\midrule
Australian & 30.00 & 0.50  \\
Phishing    & 1228.33 & 0.10  \\
Mushroom   & 8124.00 & 0.01  \\
Splice    &1000.00 & 0.90 \\
\bottomrule
\end{tabular}
\end{table}

\begin{figure}
\centering
\begin{subfigure}{.5\textwidth}
  \centering
  \includegraphics[width=.99\linewidth]{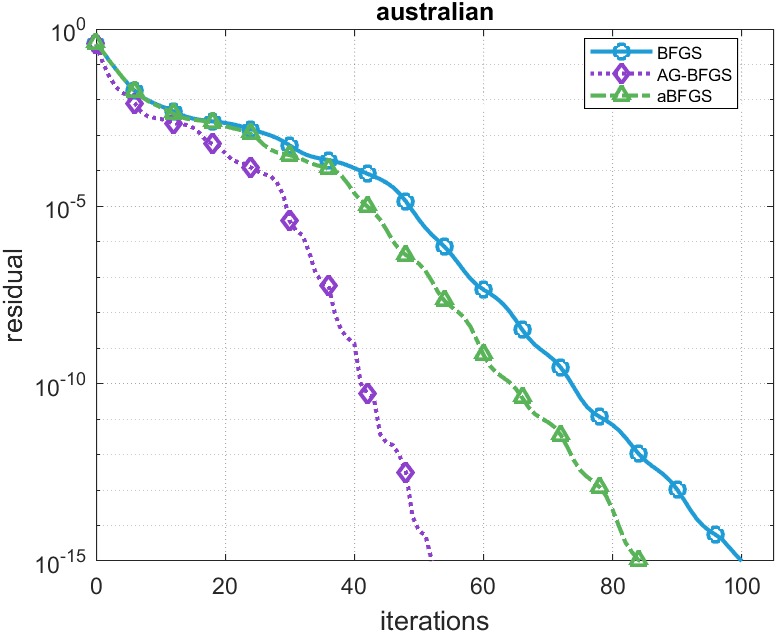}
  \caption{Australian dataset}
  \label{fig:australian}
\end{subfigure}%
\begin{subfigure}{.5\textwidth}
  \centering
  \includegraphics[width=.99\linewidth]{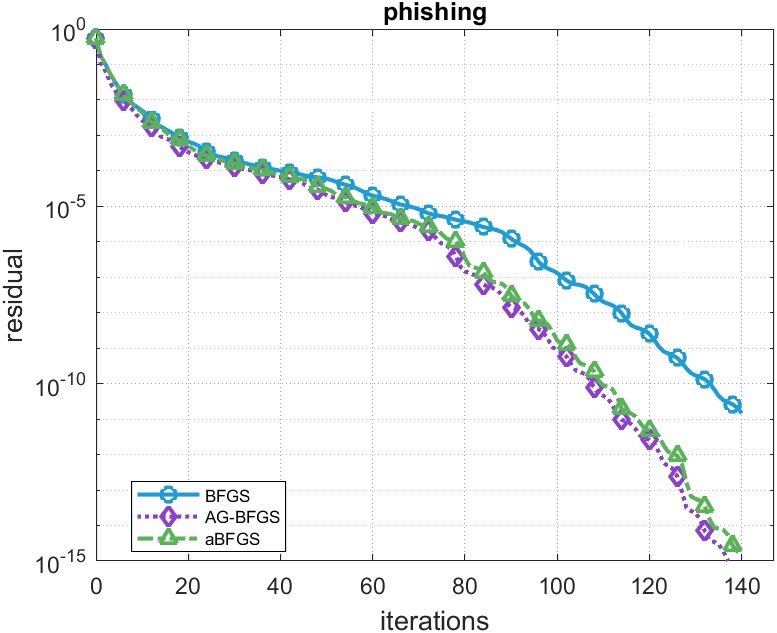}
  \caption{Phishing dataset}
  \label{fig:phishing}
\end{subfigure}
\caption{AG-BFGS vs aBFGS vs standard BFGS using the corresponding parameters shown in Table \ref{tab:parameters_values}.}
\label{fig:AP}
\end{figure}

\begin{figure}
\centering
\begin{subfigure}{.5\textwidth}
  \centering
  \includegraphics[width=.99\linewidth]{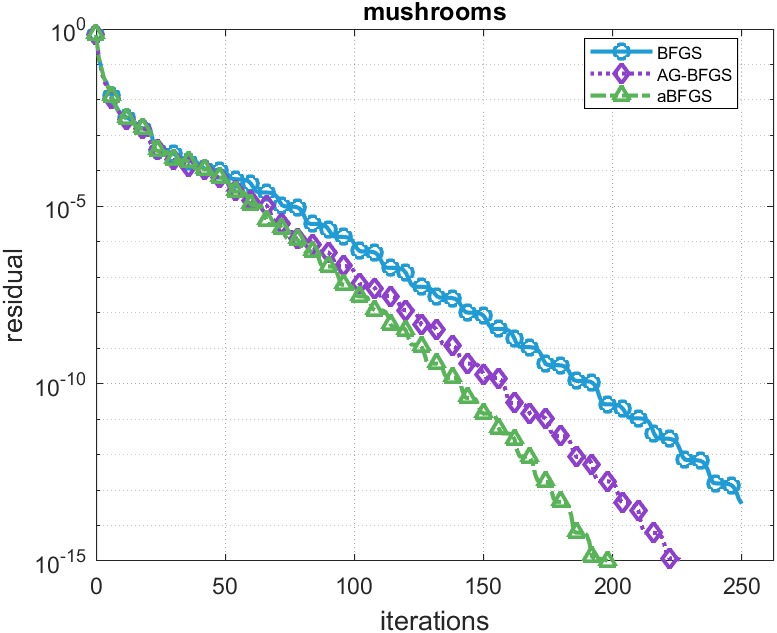}
  \caption{Mushrooms dataset}
  \label{fig:mushroom}
\end{subfigure}%
\begin{subfigure}{.5\textwidth}
  \centering
  \includegraphics[width=.99\linewidth]{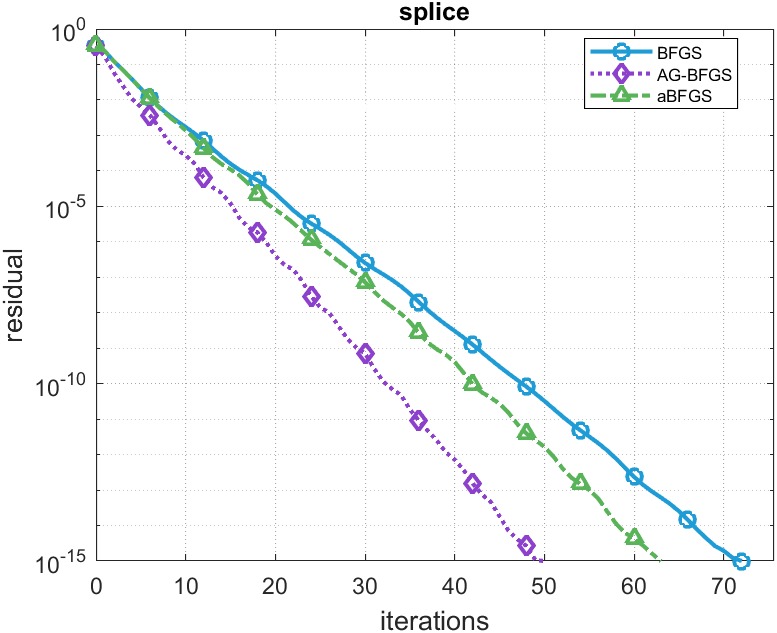}
  \caption{Splice dataset}
  \label{fig:splice}
\end{subfigure}
\caption{AG-BFGS vs aBFGS vs standard BFGS using the corresponding parameters shown in Table \ref{tab:parameters_values}.}
\label{fig:AP2}
\end{figure}

Figures~\ref{fig:australian}--\ref{fig:splice} present the experimental results for all datasets. In the figure legends, aBFGS-$\mu$-$\nu$-$\eta$ denotes the accelerated BFGS method proposed in \cite{gower2018accelerated}, parameterized by $\mu$ and $\nu$. Similarly, AG-BFGS denotes the accelerated geometry BFGS method introduced in this work, with parameters $\rho$, $\beta$, and $\gamma$ given by 0.35, 0.75, and 0.5, respectively.  The label BFGS refers to the classical BFGS algorithm, which is used as a baseline for comparison.

The results show that both accelerated variants generally outperform the standard BFGS method. When comparing aBFGS and AG-BFGS, however, it is important to recall the different parameter-selection strategies. The aBFGS parameters are problem-specific and were manually selected to achieve the best performance, whereas AG-BFGS employs the generic midpoint strategy described above, without any problem-specific information or parameter tuning. Despite this, AG-BFGS outperforms aBFGS on three of the four datasets considered, \texttt{australian}, \texttt{phishing}, and \texttt{spline}, while aBFGS performs better only on \texttt{mushrooms}. These results are particularly encouraging, as they suggest that AG-BFGS can achieve competitive or superior performance without relying on problem-dependent parameter selection.

\subsection{High-conditioning and Unimodal Problems}
\label{subsec:unimodal}

In the following, we consider five noise-free, real-parameter, single-objective benchmark functions and their rotated variants. All functions are defined on $\mathbb{R}^D$, with search domain $\mathbf{x}\in[-100,100]^D$, and have a global optimum at $\mathbf{x}^{}=\mathbf{0}$ satisfying $f(\mathbf{x}^{})=0$. The rotated versions are obtained through the orthogonal transformation $\mathbf{z}=Q^{\top}\mathbf{x}$, where $Q\in\mathbb{R}^{D\times D}$ is an orthogonal matrix satisfying $Q^{\top}Q=QQ^{\top}=I$. In all experiments, $Q$ is generated via a QR decomposition of a Gaussian random matrix and then kept fixed. This transformation preserves the location and value of the optimum while rendering the functions non-separable. Further details on these benchmark functions can be found in \cite{Hansen09}.

\vspace{10pt}

\noindent{\textbf{Ellipsoidal Function}}: Globally quadratic and ill-conditioned (condition number, $c_n$, $10^6$) function  with smooth local irregularities.

\begin{equation*}\begin{array}{l}
f_E(\mathbf{x}) = \displaystyle \sum_{i=1}^D~ 10^{6 \frac{i-1}{D-1}} x_i^2.  \\ 
\end{array}
\end{equation*}

\vspace{10pt}

\noindent{\textbf{Discuss Function}}: Globally quadratic and ill-conditioned ($c_n=10^6$) function with local irregularities. A single direction in search space is a thousand times more sensitive than others. (effect of constraint-like penalization)

\begin{equation*}\begin{array}{l}
f_D(\mathbf{x}) = 10^6x_1^2 +\displaystyle \sum_{i=2}^D~  x_i^2 .  \end{array}
\end{equation*}

\vspace{10pt}

\noindent{\textbf{Bent Cigar Function}}:  This problem has a smooth but narrow ridge along which $\sum_{i=2}^D~ x_i^2 =0$ needs to be followed. It is an ill-conditioned ($c_n=10^6$) function.

\begin{equation*}\begin{array}{l}
f_{BC}(\mathbf{x}) =  x_1^2 + 10^6\displaystyle \sum_{i=2}^D~  x_i^2.
\end{array}
\end{equation*}

\vspace{10pt}

\noindent{\textbf{Sharp Ridge Function}}: Similar to the Bent Cigar function, this problem requires following a ridge defined by $\sum_{i=2}^D z_i^2 = 0$. The ridge creates a narrow valley with a sharp transition in curvature, making the optimization process challenging, particularly for methods that rely on local curvature information.


\begin{equation*}\begin{array}{l}
f_{SR}(\mathbf{x}) =  x_1^2 + 100\displaystyle \sqrt{\sum_{i=2}^D~  x_i^2}. \end{array}\end{equation*}

\vspace{10pt}

\noindent{\textbf{Different Powers Function}}: Due to the different exponents, the sinsitivies of the $x_i$-variables become more and more different when approaching the optimum. 
As a consequence, the function exhibits heterogeneous
curvature across the search space.

\begin{equation*}\begin{array}{l}
f_{DP}(\mathbf{x}) = \displaystyle \sqrt{\sum_{i=1}^D~  |x_i|^{2+4 \frac{1-i}{D-1}}}. \end{array}
\end{equation*}

\begin{figure}
\centering
\begin{subfigure}{.5\textwidth}
  \centering
  \includegraphics[width=.9\linewidth]{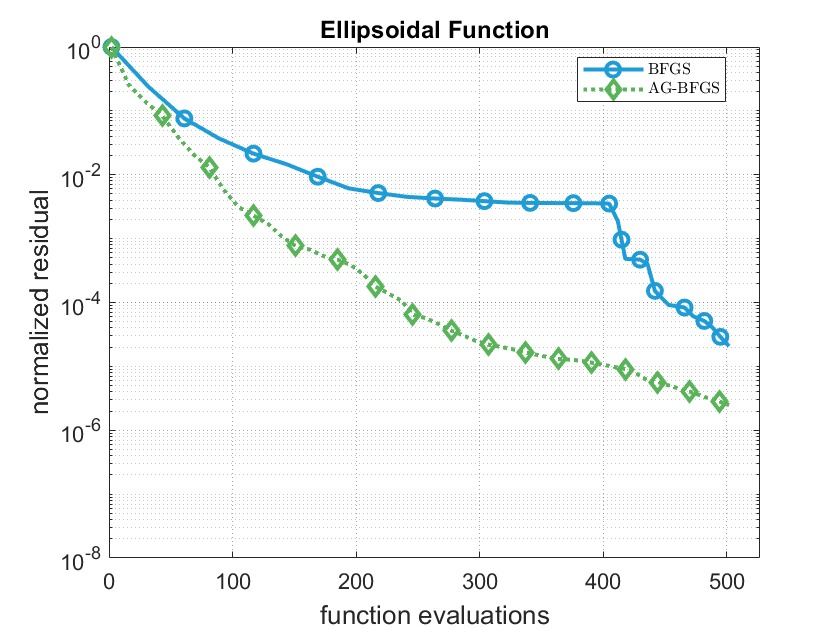}
  \caption{Elipsoidal Function}
  \label{fig:elip}
\end{subfigure}%
\begin{subfigure}{.5\textwidth}
  \centering
  \includegraphics[width=.9\linewidth]{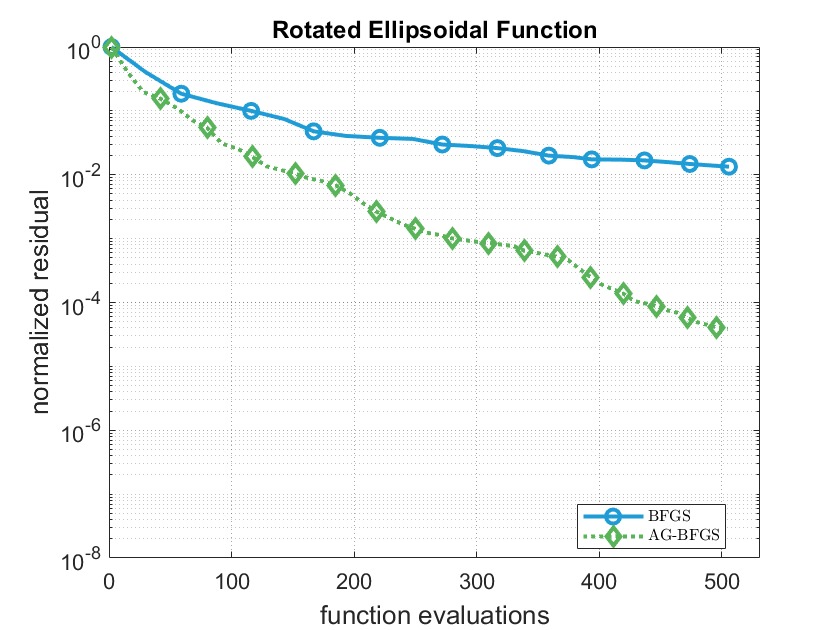}
  \caption{Rotated Ellipsoidal Function}
  \label{fig:rot_elip}
\end{subfigure}
\caption{Optimization trajectories of BFGS and the proposed accelerated BFGS (AG-BFGS) on the Ellipsoidal function. (a) shows the original separable function, while (b) shows its rotated non-separable version.}
\label{fig:test}
\end{figure}

\begin{figure}
\centering
\begin{subfigure}{.5\textwidth}
  \centering
  \includegraphics[width=.9\linewidth]{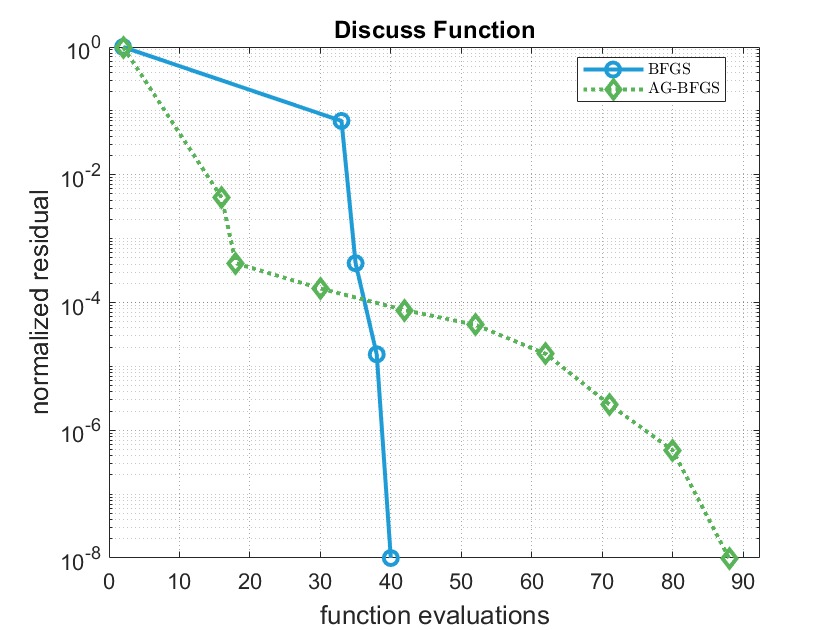}
  \caption{Discuss Function}
  \label{fig:discuss}
\end{subfigure}%
\begin{subfigure}{.5\textwidth}
  \centering
  \includegraphics[width=.9\linewidth]{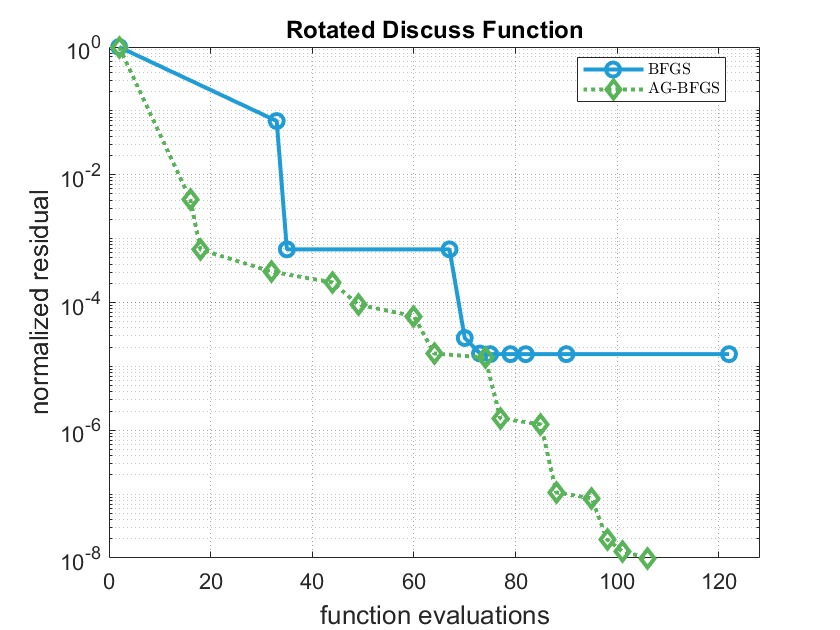}
  \caption{Rotated Discuss Function}
  \label{fig:rot_discuss}
\end{subfigure}
\caption{Optimization trajectories of BFGS and the proposed accelerated BFGS (AG-BFGS) on the Discuss function.  (a) shows the original separable function, while (b) shows its rotated non-separable version.}
\label{fig:test2}
\end{figure}

\begin{figure}
\centering
\begin{subfigure}{.5\textwidth}
  \centering
  \includegraphics[width=.9\linewidth]{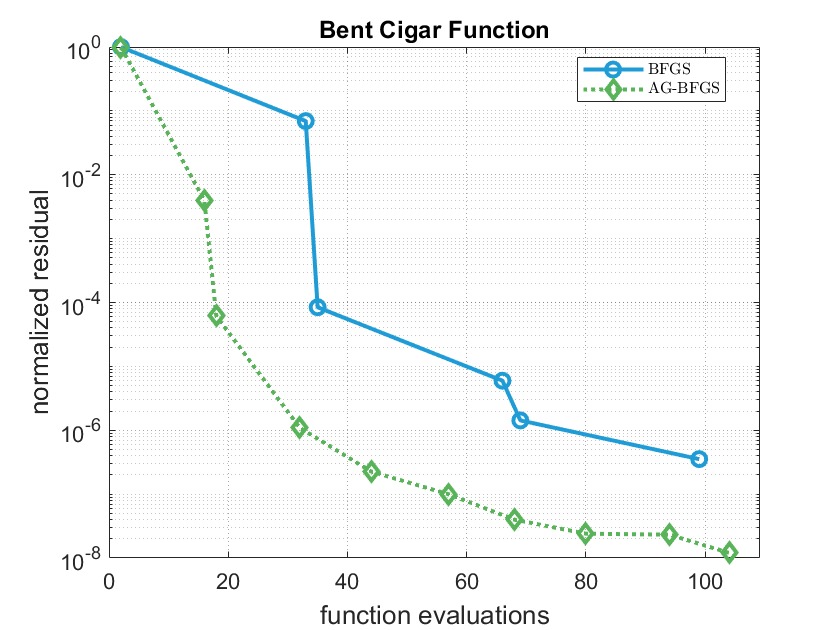}
  \caption{Bent Cigar Function}
  \label{fig:bent}
\end{subfigure}%
\begin{subfigure}{.5\textwidth}
  \centering
  \includegraphics[width=.9\linewidth]{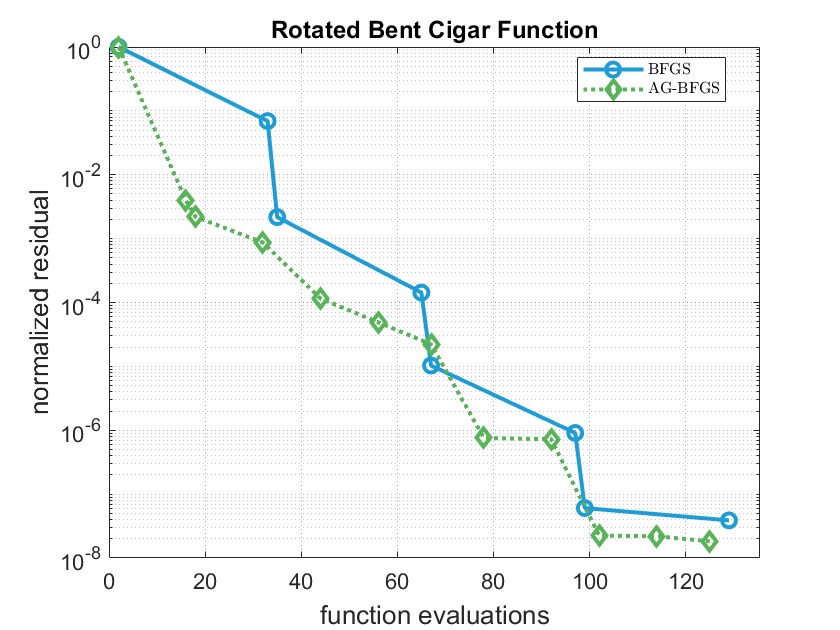}
  \caption{Rotated Bent Cigar Function}
  \label{fig:rot_bent}
\end{subfigure}
\caption{Optimization trajectories of BFGS and the proposed accelerated BFGS (AG-BFGS) on the Bent Cigar function.  (a) shows the original separable function, while (b) shows its rotated non-separable version.}
\label{fig:test3}
\end{figure}

\begin{figure}
\centering
\begin{subfigure}{.5\textwidth}
  \centering
  \includegraphics[width=.9\linewidth]{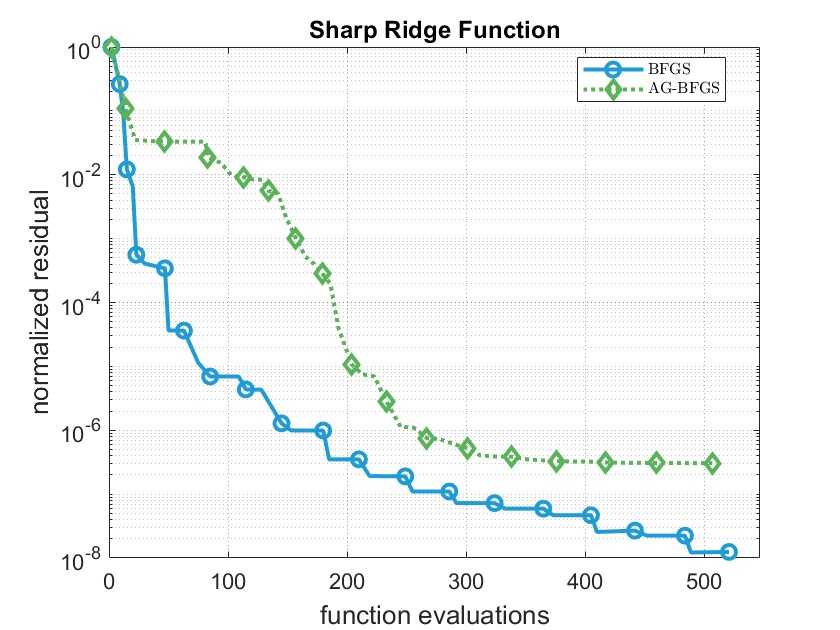}
  \caption{Sharp Ridge Function}
  \label{fig:sharp}
\end{subfigure}%
\begin{subfigure}{.5\textwidth}
  \centering
  \includegraphics[width=.9\linewidth]{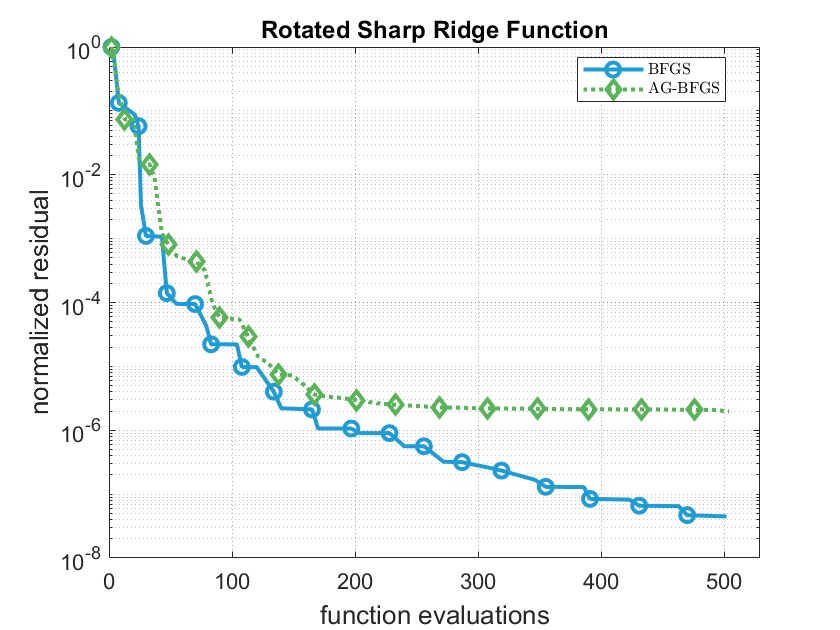}
  \caption{Rotated Sharp Ridge Function}
  \label{fig:rot_sharp}
\end{subfigure}
\caption{Optimization trajectories of BFGS and the proposed accelerated BFGS (AG-BFGS) on the Sharp Ridge function.  (a) shows the original separable function, while (b) shows its rotated non-separable version.}
\label{fig:test4}
\end{figure}

Figures \ref{fig:test} - \ref{fig:test4} show the comparison between the classical BFGS and the proposed accelerated version, AG-BFGS, using the single-objective benchmark functions. The vertical axis reports the normalized residual on a logarithmic scale, while the horizontal axis shows the number of function evaluations. For each function, the parameters, $\rho$, $\beta$ and $\gamma$, were selected using the generic mid-point strategy as described previously ($[\rho,~\beta,~\gamma] = [0.35,~0.75,~0.5]$). The stopping criteria was set as reaching a normalized residual of $1.00E-8$ or a maximum number of iterations. For each function, its separable and rotated non-separable versions are presented.

The results show that AG-BFGS provides a clear improvement over classical BFGS for several of the benchmark functions, despite using the same generic, problem-independent parameter setting, $[\rho,\beta,\gamma]=[0.35,0.75,0.5]$, in all experiments. The improvement is particularly pronounced for the Bent Cigar and Different Powers functions, for which AG-BFGS reduces the normalized residual substantially faster and reaches the stopping tolerance with fewer function evaluations. This behaviour is retained under rotation, indicating that the acceleration is not restricted to the separable structure of these problems. For the Ellipsoidal function, AG-BFGS also achieves a consistently lower residual than BFGS in both the separable and rotated cases; although neither method reaches the prescribed tolerance within the evaluation limit, the advantage of AG-BFGS is especially marked for the rotated problem. The results for the Discuss function are more nuanced: BFGS reaches the tolerance considerably faster for the separable version (approximately 40 function evaluations, compared with about 90 for AG-BFGS), whereas rotation reverses the relative behaviour. In the rotated Discuss problem, BFGS stagnates at a residual of approximately $10^{-5}$, while AG-BFGS continues to make progress and reaches the $10^{-8}$ tolerance.

The main exception is the Sharp Ridge function, where classical BFGS performs better for both the separable and rotated versions. A plausible explanation lies in the non-smooth geometry of this problem. Unlike the smooth quadratic benchmarks, the Sharp Ridge function is non-differentiable on the ridge, and curvature information becomes unreliable as the iterates approach it. Since the acceleration mechanism introduces additional dependence on information from previous iterates, it may amplify the adverse effect of inaccurate local curvature information and hinder progress along the ridge. Classical BFGS, despite also relying on curvature approximations, appears more robust in this setting. For the separable Discuss function, the opposite behaviour can be explained by its very simple but highly anisotropic structure: only one coordinate has the large curvature coefficient, and the coordinate axes are aligned with the eigendirections of the Hessian. Classical BFGS can therefore identify and adapt to this structure very efficiently, leaving comparatively little scope for acceleration and potentially making the additional momentum-like mechanism counterproductive. Once the Discuss function is rotated, this favourable coordinate alignment disappears and the curvature becomes coupled across variables.


\begin{figure}
\centering
\begin{subfigure}{.5\textwidth}
  \centering
  \includegraphics[width=.9\linewidth]{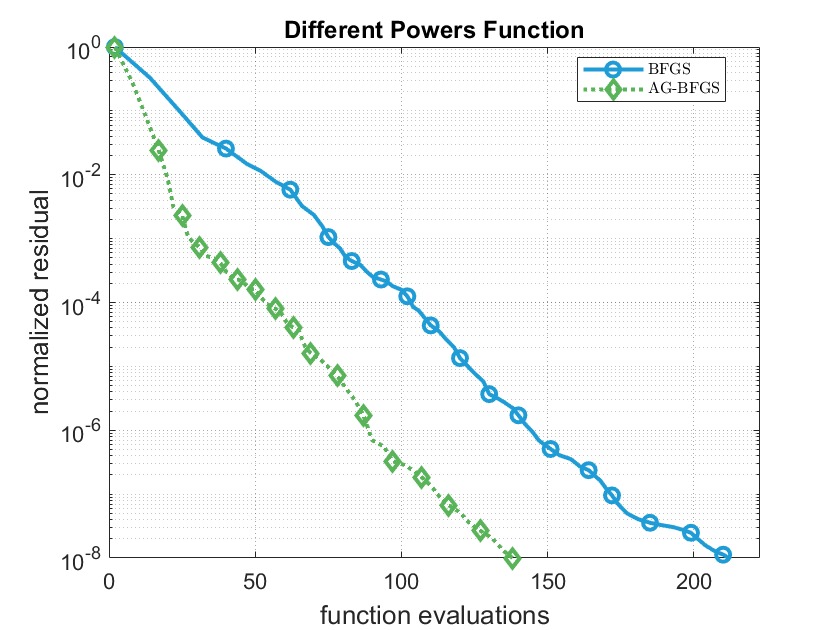}
  \caption{Different Powers Function}
  \label{fig:powers}
\end{subfigure}%
\begin{subfigure}{.5\textwidth}
  \centering
  \includegraphics[width=.9\linewidth]{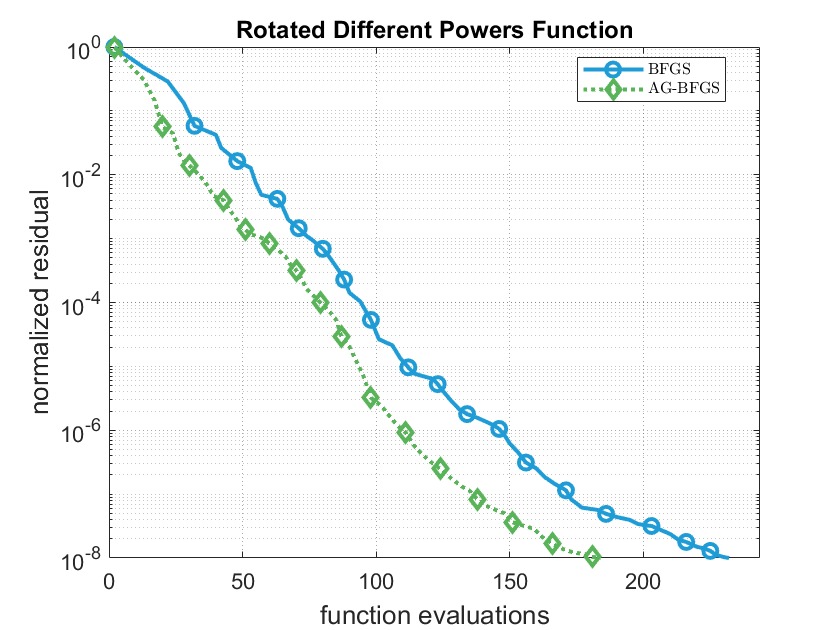}
  \caption{Rotated Different Powers Function}
  \label{fig:rot_powers}
\end{subfigure}
\caption{Optimization trajectories of BFGS and the proposed accelerated BFGS (AG-BFGS) on the Different Powers function. (a) shows the original separable function, while (b) shows its rotated non-separable version.}
\label{fig:test5}
\end{figure}

%% file: conclusion.tex
\section{Conclusion}
\label{sec:conclusion}

We have proposed an accelerated quasi-Newton method that incorporates inertial dynamics into the BFGS framework. A theoretical analysis establishes convergence to stationary points under standard assumptions, while numerical experiments demonstrate improved performance on a suite of logistic regression problems.

We further have shown that, under a suitable choice of parameters, the proposed acceleration leads to significant performance gains over the standard BFGS algorithm. This work confirms and provides theoretical support for the claim in \cite{gower2018accelerated} that Nesterov-type acceleration can be successfully applied to quasi-Newton methods.